\documentclass[12pt,twoside]{article}
\usepackage{amsmath,amssymb,amsthm}

\theoremstyle{plain}
\newtheorem{thm}{Theorem}[section]
\newtheorem{prop}{Proposition}[section]
\newtheorem{lem}{Lemma}[section]

\theoremstyle{remark}
\newtheorem{rem}{Remark}[section]

\newtheorem*{pf}{Proof}

\def\P{{\mathbb {P}}}                       
\def\E{{\mathbb {E}}}                       \def\D{{\mathbb {D}}}

\def\Z{{\mathbb {Z}}}

\def\FD{{\mathcal F}}

\def\o{{\rm {o}}} 
\def\O{{\rm {O}}} 

\def\CC{{\rm\kern.24em
   \vrule width.02em
       height1.4ex depth-.05ex
   \kern-.26em C}}
\def\QQ{{\rm\kern.24em
   \vrule width.02em
       height1.4ex depth-.05ex
   \kern-.26em Q}}
\def\PP{{\rm I\kern-.25em P}}         \def\RR{{\rm I\kern-.25em R}}
\def\DD{{\rm I\kern-.25em D}}         \def\EE{{\rm I\kern-.25em E}}
\def\FF{{\rm I\kern-.25em F}}         \def\NN{{\rm I\kern-.23em N}}
\def\RRp{{\rm I\kern-.25em R}_{+}}
\def\IND{{\rm 1\kern-.25em I}}

\begin{document}

\centerline{\large{\bf {Scale-free property for degrees and weights}}}
\smallskip
\centerline{\large{\bf {in an $N$-interactions random graph model}}}

\bigskip

\centerline{ {\sc Istv\'an Fazekas} and {\sc Bettina Porv\'azsnyik} }

\bigskip
Faculty of Informatics, University of Debrecen,  P.O. Box 12, 4010 Debrecen, Hungary, 
\centerline{
e-mail: fazekas.istvan@inf.unideb.hu, porvazsnyik.bettina@inf.unideb.hu.}

\medskip \bigskip

\begin{abstract}
A general random graph evolution mechanism is defined.
The evolution is a combination of the preferential attachment model and the interaction of $N$ vertices $\left( N \geq 3\right)$.
A vertex in the graph is characterized by its degree and its weight.
The weight of a given vertex is the number of the interactions of the vertex.
The asymptotic behaviour of the graph is studied.
Scale-free properties both for the degrees and the weights are proved.
It turns out that any exponent in $(2,\infty)$ can be achieved.
The proofs are based on discrete time martingale theory.
\end{abstract}

\renewcommand{\thefootnote}{}
{\footnotetext{
{\bf Key words and phrases:}
Random graph, preferential attachment, scale-free, power law, submartingale, Doob-Meyer decomposition.

{\bf Mathematics Subject Classification:} 05C80, 
60G42. 

The publication was supported by the T\'AMOP-4.2.2.C-11/1/KONV-2012-0001 project. The project has been supported by the European Union, co-financed by the European Social Fund.
}

\section{Introduction}
Network theory became a popular field during the last $15$ years.
Several real-world networks were investigated such as the WWW, the Internet, social and biological networks (see \cite{durrett} for an overview).
It turned out that a main common characteristic of such networks is their scale-free nature, in other words the asymptotic power law degree distribution, 
i.e.  $p_k \sim Ck^{-\gamma}$, as $k\to\infty$.
Using large data sets, it was shown that for the WWW the in-degree and the out-degree of web pages follow power law with $\gamma_{\rm{in}}=2.1$ and $\gamma_{\rm{out}}=2.7$, 
for the Internet $\gamma=2.3$, for the movie actor network $\gamma=2.3$, for the collaboration graph of mathematicians $\gamma=2.4$ (see \cite{durrett} for details).
To describe the phenomenon, in \cite{barabasi} the preferential attachment model was suggested.
In the preferential attachment model the growing mechanism of the random graph is the following. 
At every time $t=2,3,\dots$ a new vertex with $m$ edges is added so that the edges link the new vertex to $m$ old vertices.
The probability $\pi_i$ that the new vertex will be connected to the old vertex $i$ depends on the degree $d_i$ of vertex $i$, so that $\pi_i= d_i/\sum_j d_j$.
The power law degree distribution in the preferential attachment model was proved by a couple of methods (see, e.g. \cite{bollobas}).

There are several modifications of the preferential attachment model (see \cite{cooper}, \cite{SGWN}).
It is also known that besides the degrees of the vertices other characteristics of the graph can be important (see \cite{SGWN}).
In \cite{BaMo1} a model based on the interaction of three vertices was introduced.
The power law degree distribution in that model was proved in \cite{BaMo2}.
In \cite{FIPB}, instead of the three-interactions model, interactions of four vertices were studied.
It turned out that in the seemingly complicated four-interactions model the asymptotic behaviour is as simple as in the three-interactions model.
Therefore it is hopeful that the overburdening formulae of the $N$-interactions model lead to tractable asymptotic results.

In this paper, we extend the model and the results of \cite{BaMo1}, \cite{BaMo2} and \cite{FIPB} to interactions of $N$ vertices.
Our model is the following.
A complete graph with $m$ vertices we call an $m$-clique, for short.
We denote an $m$-clique by the symbol $K_m$.
In our model at time $n=0$ we start with a $K_N$.
The initial weight of this graph is one. 
This graph contains $N$ vertices, $ \binom{N}{2}$ edges, \dots , $\binom{N}{M}$ $M$-cliques $\left( M \leq N \right)$. 
Each of these objects has initial weight $1$.
After the initial step we start to increase the size of the graph. 
At each step, the evolution of the graph is based on the interaction of $N$ vertices. 
More precisely, at each step $n=1,2,\dots$ we consider $N$ vertices and draw all non-existing edges between these vertices.
So we obtain a $K_N$.
The weight of this graph $K_N$ and the weights of all cliques in $K_N$ are increased by $1$. 
(That is we increase the weights of $N$ vertices, $ \binom{N}{2}$ edges, \dots , $N$ different $\left(N-1\right)$-cliques and the $N$-clique $K_N$ itself.)
The choice of the $N$ interacting vertices is the following.

There are two possibilities at each step.
With probability $p$ we add a new vertex that interacts with $N-1$ old vertices, on the other hand,
with probability $\left( 1-p \right)$, $N$ old vertices interact. 
Here $0 < p \leq 1$ is fixed.

When we add a new vertex, then we choose $N-1$ old vertices and they together will form an $N$-clique.
However, to choose the $N-1$ old vertices we have two possibilities.
With probability $r$ we choose an $\left(N-1\right)$-clique from the existing $\left(N-1\right)$-cliques according to the weights of the $\left(N-1\right)$-cliques.
It means that an $\left(N-1\right)$-clique of weight $w_t$ is chosen with probability $w_t/\sum_h w_h$.
On the other hand, with probability $1-r$, we choose among the existing vertices uniformly, that is all $N-1$ vertices have the same chance.

At a step when we do not add a new vertex, then $N$ old vertices interact. 
As in the previous case, we have two possibilities. 
With probability $q$, we choose one $K_N$ of the existing $N$-cliques according to their weights.
It means that an $N$-clique of weight $w_t$ is chosen with probability $w_t/\sum_h w_h$.
On the other hand, with probability $1-q$, we choose among the existing vertices uniformly, that is all subsets consisting of $N$ vertices have the same chance.

In this paper we show that the above mechanism results in a scale-free graph.
To prove our results we follow the lines of \cite{BaMo1}, \cite{BaMo2}.
Let $X(n,d,w)$ denote the number of vertices of weight $w$ and degree $d$ after the $n$th step.
Let $V_n$ denote the number of vertices after the $n$th step.
Let $\FD_{n-1}$ denote the $\sigma$-algebra of observable events after the $(n-1)$th step.
First we calculate the conditional expectation ${\E} \{ X(n,d,w)|\FD_{n-1} \}$, see Lemma \ref{EX/F}.
Then we prove (Theorem \ref{limX/V}) that the ratio $\frac{X \left( n,d,w \right)}{V_n}$ converges to $x_{d,w}$
almost surely (a.s.) as $n \rightarrow \infty $, where the limits $x_{d,w}$ are fixed non-negative numbers.
The main tool of the proof is the Doob-Meyer decomposition of submartingales.
We remark that in the $3$-interactions model the limits $x_{d,w}$ are always positive (see \cite{BaMo2}).
However, in the $N$-interactions model the limits $x_{d,w}$ can be zero unless $N$ is equal to $3$.
It is an important phenomenon, because the appearance of zero limits simplifies the seemingly intractable formulae.

We show that $x_{d,w}$, $d=N-1,N,\dots, \left(N-1\right)w$, $w=1,2,\dots$, is a proper two-dimensional discrete probability distribution (Lemma \ref{xdw}).
Then we turn to the scale-free property for the weights.
Let $X \left( n,w \right)$ denote the number of vertices of weight $w$ after the $n$th step.
Then for all $w=1,2,\dots$ we have
$$
\dfrac{X \left( n,w \right)}{V_n} \rightarrow x_{w} = x_{N-1,w} + x_{N,w} + \dots + x_{\left(N-1\right)w,w}
$$
almost surely.
Moreover, $x_{w} \sim C  w^{- \left( 1 + \frac{1}{\alpha} \right)}$, as $w\to\infty$
(Theorem \ref{theorem:scalefreeWeights}),
that is the distribution $x_{w}$ has a tail which decays as a power-law with exponent $ 1 + \frac{1}{\alpha}$.
To derive the above results from Theorem \ref{limX/V}, we need only some known facts about the $\varGamma$-function, 
see the proofs of Lemma \ref{xdw} and Theorem \ref{theorem:scalefreeWeights}.
Finally, we obtain the scale-free property for the degrees. 
Let us denote by $U\left( n,d \right)$ the number of vertices of degree $d$ after the $n$th step. 
For any $d \geq N-1$ we have 
\begin{equation*}
\dfrac{U\left( n,d \right)}{V_n} \to u_d
\end{equation*}
a.s. as $n \to \infty$, where $u_d$, $d=N-1,N, N+1,\dots$, are positive numbers. 
Furthermore,
\begin{equation*}
u_d \sim \dfrac{\varGamma \left( 1 + \frac{\beta + 1}{\alpha} \right)}{\alpha_2 \varGamma 
\left( 1 + \frac{\beta}{\alpha} \right)}\left( \dfrac{\alpha  d}{\alpha_2} \right)^{- \left( 1 + \frac{1}{\alpha} \right)},
\end{equation*}
as $d \to \infty$, where $\alpha$, $\beta$ and $\alpha_2$ are appropriate constants (see Theorem \ref{ThmScaleFreeDegree}).
We see that in both cases the exponent is $1 + \frac{1}{\alpha}$.
\section{The evolution of the graph}
Throughout the paper $0<p \leq 1$, $0 \leq r \leq 1$, $0 \leq q \leq 1$ are fixed numbers.
Let $X(n,d,w)$ denote the number of vertices of weight $w$ and degree $d$ after the $n$th step.
Let $V_n$ denote the number of vertices after the $n$th step.
\begin{rem}
Each vertex has initial weight $1$ and initial degree $N-1$. 
When a vertex takes part in an interaction, the weight of this vertex is increased by $1$ and the degree of this vertex may increase by $0, 1, 2, \dots $ or $N-1$. 
Therefore $X(n,d,w)$ can be positive only for $1 \leq w \leq n+1$ and $N-1 \leq d \leq \left(N-1\right)w$. 
\end{rem}
Let $\FD_{n-1}$ denote the $\sigma$-algebra of observable events after the $(n-1)$th step.
We compute the conditional expectation of $X(n,d,w)$ with respect to $\FD_{n-1}$ for $w \geq 1$.
The results of this paper will be based on it. 
The particular cases $N=3$ and $N=4$ are included in \cite{BaMo2} and \cite{FIPB}, respectively.

Let
\begin{equation*}
\alpha_1 = \left(1-p\right) q,  \quad \alpha_2 = \dfrac{N-1}{N}pr, \quad \alpha = \alpha_1 + \alpha_2, \quad
\beta =  \left( N-1 \right)\left( 1-r \right) + \dfrac{N\left( 1-p \right)\left( 1-q \right)}{p}.
\end{equation*} 
\begin{lem} \label{EX/F}
One has
\begin{equation*}
{\E} \{ X(n,d,w)|\FD_{n-1} \} =
X(n-1,d,w) \left[ 1-\left(\dfrac{w}{n}\alpha + \dfrac{p}{V_{n-1}}\beta \right) \right] +
\end{equation*}
\begin{equation*}
+X(n-1,d,w-1) \left[ \left( 1-p \right) \left( q\dfrac{w-1}{n} + \left( 1-q \right)\dfrac {\binom{d}{N-1}} {\binom{V_{n-1}}{N}} \right) \right] +
\end{equation*}
\begin{equation*}
+ X(n-1,d-1,w-1) \left[ p \left(r \dfrac{\left(N-1\right)\left( w-1 \right)}{Nn} + \left( 1-r \right) \dfrac{\binom{d-1}{N-2}}{\binom{V_{n-1}}{N-1}}  \right) + \right.
\end{equation*}
\begin{equation*}
\left.
+ \left( 1-p \right)\left( 1-q \right)\dfrac {\binom{d-1}{N-2} \left( V_{n-1}-d \right)} {\binom{V_{n-1}}{N}} \right] + \dots +
\end{equation*}
\begin{equation*}
+ X(n-1,d-m,w-1) \left[ p \left( 1-r \right)\dfrac { \binom{d-m}{N-m-1} \binom{V_{n-1}-d+m-1}{m-1}} {\binom{V_{n-1}}{N-1}} + \right.
\end{equation*}
\begin{equation*}
\left.
+ \left( 1-p \right)\left( 1-q \right)\dfrac {\binom{d-m}{N-m-1} \binom{ V_{n-1}-d+m-1}{m} } {\binom{V_{n-1}}{N}} \right] + \dots +
\end{equation*}
\begin{equation*}
+X(n-1,d-\left(N-1\right),w-1) \left[ p\left( 1-r \right) \dfrac{\binom{V_{n-1}-d+N-2}{N-2}}{\binom{V_{n-1}}{N-1}} \right.
\end{equation*}
\begin{equation} \label{CondExp}
\left. +
\left( 1-p \right) \left( 1-q \right)\dfrac{\binom{V_{n-1}-d+N-2}{N-1}}{\binom{V_{n-1}}{N}} \right] + p\delta_{d,N-1}\delta_{w,1}
\end{equation}
for $w\ge 1$ and $N-1\le d \le \left( N-1 \right)w$, $1 < m < N-1$.
Here $\delta_{k,l}$ denotes the Dirac-delta.
\end{lem}
\begin{pf}
The total weight of $N$-cliques after $(n-1)$ steps is $n$.
The total weight of $\left(N-1\right)$-cliques after $(n-1)$ steps is $Nn$.
The total weight of $\left(N-1\right)$-cliques having a fixed common vertex of weight $w$ is $\left( N-1 \right)w$.
Moreover, after $(n-1)$ steps, we have the following.
When we choose $\left(N-1\right)$ vertices randomly, then the probability that a given vertex is chosen is 
$$
\frac{\binom{V_{n-1}-1}{N-2}}{\binom{V_{n-1}}{N-1}} = \frac{N-1}{V_{n-1}}.
$$
When we choose $N$ vertices randomly, then the probability that a given vertex is chosen is 
$$
\frac{\binom{V_{n-1}-1}{N-1}}{\binom{V_{n-1}}{N}} = \frac{N}{V_{n-1}}.
$$
Therefore the probability that an old vertex of weight $w$ takes part in the interaction at step $n$ is
$$
p \left(r \dfrac{\left(N-1\right)w}{Nn} + (1-r) \dfrac{N-1}{V_{n-1}}\right) +
\left(1-p\right) \left(q \dfrac{w}{n} + (1-q) \dfrac{N}{V_{n-1}}\right) = 
\dfrac{w}{n}\alpha + \dfrac{p}{V_{n-1}}\beta, 
$$
where
\begin{equation*}
\alpha = (1-p)q + \dfrac{\left(N-1\right)pr}{N}, \quad \beta = \dfrac{1}{p} \left\{\left(N-1\right)p(1-r) + N(1-p)(1-q)  \right\}. 
\end{equation*}
A new vertex always takes part in the interaction. 
At each step with probability $p$ a new vertex with weight $1$ and with degree $\left(N-1\right)$ is born.

Consider a fixed vertex with weight $w$ and degree $d$. 
The probability that in the $n$th step 
\begin{itemize}
 \item neither its degree $d$ nor its weight $w$ change is
 $$
1-\left(\dfrac{w}{n}\alpha + \dfrac{p}{V_{n-1}}\beta \right)\,;
 $$
 \item its degree does not change but its weight is increased by 1 is
 $$
\left( 1-p \right) \left( q\dfrac{w-1}{n} + \left( 1-q \right)\dfrac {\binom{d}{N-1}} {\binom{V_{n-1}}{N}} \right)\,;
 $$
 \item both its degree and its weight are increased by 1 is
 $$
p \left(r \dfrac{\left(N-1\right)\left( w-1 \right)}{Nn} + \left( 1-r \right) \dfrac{\binom{d-1}{N-2}}{\binom{V_{n-1}}{N-1}}  \right) +
\left( 1-p \right)\left( 1-q \right)\dfrac {\binom{d-1}{N-2} \left( V_{n-1}-d \right)} {\binom{V_{n-1}}{N}}\,;
 $$
 \item its degree is increased by $m$ ($1<m<N-1$) and its weight is increased by 1 is
 $$
p \left( 1-r \right)\dfrac { \binom{d-m}{N-m-1} \binom{V_{n-1}-d+m-1}{m-1}} {\binom{V_{n-1}}{N-1}} + \left( 1-p \right)\left( 1-q \right)\dfrac {\binom{d-m}{N-m-1} \binom{ V_{n-1}-d+m-1}{m} } {\binom{V_{n-1}}{N}}\,;
 $$
 \item its degree is increased by $N-1$ and its weight is increased by 1 is
 $$
p \left( 1-r \right) \dfrac{\binom{V_{n-1}-d+N-2}{N-2}}{\binom{V_{n-1}}{N-1}} +
\left( 1-p \right) \left( 1-q \right)\dfrac{\binom{V_{n-1}-d+N-2}{N-1}}{\binom{V_{n-1}}{N}}\,.
 $$
\end{itemize}
Using the above formulae, we obtain equation  \eqref{CondExp}.
\hfill$\Box$
\end{pf}

The following theorem is an extension of \textit{Theorem 3.1} of \cite{BaMo2}, see also \textit{Theorem 2.1} of \cite{FIPB}.
For $N>3$ we shall see, that several terms are asymptotically negligible, therefore the final expressions are as simple as in the case of $N=3$.

\begin{thm} \label{limX/V}
Let $0<p<1$, $q>0$, $r>0$ and $(1-r)(1-q)>0$.
For any fixed $w$ and $d$ with $1 \leq w$ and $N-1 \leq d \leq w\left(N-1 \right)$ we have
\begin{equation}   \label{X/Vtox}
\dfrac{X \left( n,d,w \right)}{V_n} \rightarrow x_{d,w}
\end{equation}
almost surely as $n \rightarrow \infty $, where $x_{d,w}$ are fixed non-negative numbers.
Furthermore, the numbers $x_{d,w}$ satisfy the following recurrence:
$$
x_{N-1,1} = \dfrac{1}{\alpha + \beta +1} > 0, \quad \quad \quad x_{d,1}=0, \, \text{ for } \, d\ne N-1,
$$
\begin{equation} \label{rekurzio_x(d,w)}
x_{d,w} = \dfrac{1}{\alpha w + \beta +1} \left[ \alpha_1 \left( w-1 \right) x_{d,w-1} + \alpha_2 \left(  w-1\right)x_{d-1,w-1} 
+\beta x_{d-\left( N-1 \right),w-1} \right],
\end{equation}
for $w\ge 2$, $N-1\le d\le w\left(N-1 \right)$, where
\begin{equation*}
\alpha_1 = \left(1-p\right) q, \quad \alpha_2 = \dfrac{N-1}{N}pr, \quad \alpha = \alpha_1 + \alpha_2, \quad \beta =  \left( N-1 \right)\left( 1-r \right) + \dfrac{N\left( 1-p \right)\left( 1-q \right)}{p}.
\end{equation*}

If $w \geq 1$ is fixed then there exists $d$ with $N-1\le d \le w\left(N-1\right)$ such that $x_{d,w}$ is positive and if $w \geq 1$ and $N \geq 4 $ then there exists $d$ with $N-1\le d \le w\left(N-1\right)$ such that $x_{d,w}$ is equal to zero.
Moreover, in the cases when $x_{d,w} = 0$ we have
\begin{equation*}
\dfrac{X \left( n,d,w \right)}{V_n} = \o \left( n^{-a} \right),
\end{equation*}
where $a$ is a positive number which may depend on $w$ and $d$.


If $N-1\le d \le w\left(N-1\right)$ does not satisfied, then $x_{d,w}=0$.
\end{thm}

\begin{pf}
We follow the lines of \cite{BaMo2}.
Introduce notation
\begin{equation} \label{c_def}
c(n,w) = \prod_{i=w-1}^{n} \left( 1-\dfrac{\alpha w}{i}-\dfrac{\beta p}{V_{i-1}} \right)^{-1}, \quad n \geq w-1, w \geq 1.
\end{equation}
$c(n,w)$ is an $\FD_{n-1}$-measurable random variable.
Applying the Marcinkiewicz strong law of large numbers to the number of vertices, we have
\begin{equation}  \label{Markinkiewicz}
V_n = pn + \o \left( n^{1/2 + \varepsilon} \right)
\end{equation}
almost surely, for any $\varepsilon >0$.

Using \eqref{Markinkiewicz} and the Taylor expansion for $\log(1+x)$, we obtain
\begin{equation*}
\log c\left( n,w \right) = -\sum_{i=w-1}^{n} \log  \left( 1-\dfrac{\alpha w}{i}-\dfrac{\beta }{i + \o \left( i^{1/2 + \varepsilon} \right)} \right) = 
\left( \alpha w + \beta \right) \sum_{i=w-1}^{n} \dfrac{1}{i} + \O \left( 1 \right),
\end{equation*}
where the error term is convergent as $n \to \infty$.
Therefore
\begin{equation} \label{c(n,w)_asz.}
c(n,w) \sim a_{w} n^{\alpha w + \beta}
\end{equation}
almost surely, as $n \to \infty$, where $a_{w}$ is a positive random variable.

Let
\begin{equation*} \label{Z_def}
Z \left( n,d,w \right) = c\left( n,w \right) X \left( n,d,w \right) \quad \text{for} \quad 1 \leq w, \, N-1 \leq d \leq w\left(N-1 \right).
\end{equation*}
Using \eqref{CondExp}, we can see that $ \left\{ Z \left( n,d,w \right) , \FD_{n} , n=w-1,w,w+1,\dots \right\}$ is a non-negative submartingale for any fixed 
$1 \leq w$, $N-1 \leq d \leq \left( N-1 \right)w$.
Define $ Z \left( n,d,w \right) = 0$ for $n=1,2,\dots,w-2$.
Applying the Doob-Meyer decomposition to $ Z \left( n,d,w \right)$, we can write
\begin{equation*} \label{Doob_Meyer}
Z \left( n,d,w \right) = M \left( n,d,w \right) + A \left( n,d,w \right),
\end{equation*}
where $M \left( n,d,w \right)$ is a martingale and $A \left( n,d,w \right)$ is a predictable increasing process. 
The general form of $M \left( n,d,w \right)$ and $A \left( n,d,w \right)$ are the following:
\begin{equation} \label{M_alt}
M \left( n,d,w \right) = \sum_{i=1}^{n} \left[ Z \left( i,d,w \right) - {\E} \left( Z \left( i,d,w \right) | \FD_{i-1} \right) \right],
\end{equation}
\begin{equation}  \label{A_alt}
A \left( n,d,w \right) = {\E} Z \left( 1,d,w \right) + \sum_{i=2}^{n} \left[ {\E} \left( \left( Z \left( i,d,w \right) | \FD_{i-1} \right) - Z \left( i-1,d,w \right) \right) \right],
\end{equation}
where $\FD_{0}$ is the trivial $\sigma$-algebra.
Using \eqref{A_alt} and \eqref{CondExp}, we have
\begin{multline} \label{A(n,d,w)}
\begin{split}
A & \left( n,d,w \right) = 
{\E} Z \left( 1,d,w \right)+\\
& +\sum_{i=2}^{n} \left[ c \left( i,w \right) X \left( i-1,d,w-1 \right) \left( 1-p \right) \left( q \dfrac{w-1}{i} + \left( 1-q \right) \dfrac{\binom{d}{N-1}}{\binom{V_{i-1}}{N}} \right) \right.+ \\
& +c \left( i,w \right) X(i-1,d-1,w-1) \times \\ 
&  \times \left( p \left(r \dfrac{\left(N-1\right)\left( w-1 \right)}{Ni} + \left( 1-r \right) \dfrac{\binom{d-1}{N-2}}{\binom{V_{i-1}}{N-1}}  \right) \right. +
 \left( 1-p \right) \left( 1-q \right) \left. \dfrac{\binom{d-1}{N-2} \left( V_{i-1}-d \right)}{\binom{V_{i-1}}{N}} \right) +\\
& + \dots +c\left( i,w \right)X(i-1,d-m,w-1) \times \\ 
& \times \left( p \left( 1-r \right)\dfrac { \binom{d-m}{N-m-1} \binom{V_{i-1}-d+m-1}{m-1} } {\binom{V_{i-1}}{N-1}} \right.
 + \left. \left( 1-p \right)\left( 1-q \right)\dfrac {\binom{d-m}{N-m-1} \binom{ V_{i-1}-d+m-1}{m} } {\binom{V_{i-1}}{N}} \right) +\\
& + \dots  +c\left( i,w \right) X(i-1,d-\left(N-1\right),w-1) \times \\ 
& \times \left( p \left( 1-r \right) \dfrac{\binom{V_{i-1}-d+N-2}{N-2}}{\binom{V_{i-1}}{N-1}} +
\left( 1-p \right) \left( 1-q \right) \left. \dfrac {\binom{V_{i-1}-d+N-2}{N-1}} {\binom{V_{i-1}}{N}} \right) +  c \left( i,w \right) p\delta_{d,N-1}\delta_{w,1} \right].
\end{split}
\end{multline}
Let $B \left( n,d,w \right)$ the sum of the conditional variances of $Z \left( n,d,w \right)$.
Now we give an upper bound of $B \left( n,d,w \right)$.
\begin{equation*} 
B \left( n,d,w \right) = \sum_{i=2}^{n} {\D}^2 \left( Z \left( i,d,w \right) | \FD_{i-1} \right) =
 \sum_{i=2}^{n} {\E} \{ \left( Z \left( i,d,w \right) - {\E} \left( Z \left( i,d,w \right) |\FD_{i-1}\right) \right)^2 | \FD_{i-1}  \} = 
\end{equation*}
\begin{equation*}
=
 \sum_{i=2}^{n} c\left( i,w \right)^2 {\E} \{ \left( X \left( i,d,w \right) - {\E} \left( X \left( i,d,w \right) |\FD_{i-1}\right) \right)^2 | \FD_{i-1} 
 \} \leq 
\end{equation*}
\begin{equation*}
\leq 
\sum_{i=2}^{n} c\left( i,w \right)^2 {\E} \{ \left( X \left( i,d,w \right) -  X \left( i-1,d,w \right) \right)^2 | \FD_{i-1} 
 \} \leq
\end{equation*}
\begin{equation}\label{Bndw}
\leq 
 N^2 \sum_{i=2}^{n} c\left( i,w \right)^2 = \O\left( n^{2 \left( \alpha w+ \beta \right) +1} \right).
\end{equation}
Above we used that $c\left( i,w \right)$ is $\FD_{i-1}$-measurable and at each step $N$ vertices can interact.
Finally, we applied \eqref{c(n,w)_asz.}.
Jensen's inequality implies that $M^2 \left( n,d,w \right)$ is a (non-negative) submartingale if $M \left( n,d,w \right)$ is a martingale. 
Now we can apply the Doob-Meyer decomposition to $M^2 \left( n,d,w \right)$.
It is known that $B(n, d,w)$, that is the sum of the conditional variances of terms $Z(n,d,w)$ from formula \eqref{Bndw}, 
is the same (up to an additive constant) as the increasing predictable process in the Doob-Meyer decomposition of the non-negative submartingale $M^2(n,d,w)$.
Therefore the Doob-Meyer decomposition is
$$
M^2 \left( n,d,w \right) = Y \left( n,d,w \right) + B \left( n,d,w \right),
$$
where $Y \left( n,d,w \right)$ is a martingale and the predictable increasing process $B \left( n,d,w \right)$ is given by \eqref{Bndw}.

We use induction on $w$.
Let $w=1$. 
We can see that a vertex of weight $1$ could take part in an interaction when it was born. 
Therefore its degree must be equal to $N-1$. 
By \eqref{A(n,d,w)},
\begin{equation} \label{Teljes indukcio:w=1}
A \left( n,N-1,1 \right) \sim p \sum_{i=2}^{n} c\left( i,1 \right) \sim p \sum_{i=2}^{n} a_{1} i^{\alpha + \beta} \sim p a_{1} \dfrac{n^{\alpha + \beta + 1}}{\alpha + \beta + 1}
\end{equation}
a.s. as $n \to \infty$. \\
By \eqref{Bndw}, $B \left( n,N-1,1 \right) = \O \left( n^{2\left( \alpha + \beta \right) + 1} \right)$ and therefore 
$ \left(B \left( n,N-1,1 \right)\right)^{\frac{1}{2}} \log B \left( n,N-1,1 \right) = \,\,\,   = \O \left( A \left( n,N-1,1 \right)\right)$.
Therefore, by Proposition VII-2-4 of \cite{neveu},
\begin{equation} \label{Z(n,3,1)_A}
Z \left( n,N-1,1 \right) \sim A \left( n,N-1,1 \right) \quad \text{a.s. on the event } \quad \{A \left( n,N-1,1 \right) \to \infty\} \quad \text{as}\quad n \to \infty.
\end{equation}
As, by \eqref{Teljes indukcio:w=1}, $A(n,N-1,1) \to \infty$ a.s., therefore
using \eqref{Teljes indukcio:w=1}, \eqref{Markinkiewicz} and \eqref{c(n,w)_asz.}, relation \eqref{Z(n,3,1)_A} implies
\begin{equation}
\dfrac{X \left( n,N-1,1 \right)}{V_n} = \dfrac{Z \left( n,N-1,1 \right)}{c\left( n,1 \right)V_n} \sim \dfrac{A \left( n,N-1,1 \right)}{c\left( n,1 \right)V_n} 
\sim \dfrac{p a_{1} \dfrac{n^{\alpha + \beta +1}}{\alpha +\beta +1}}{a_{1} n^{\alpha + \beta }p n} = \dfrac{1}{\alpha + \beta +1} = x_{N-1,1} > 0
\end{equation}
almost surely.
So \eqref{X/Vtox} is valid for $w=1$.

Now let $w=2$. 
In this case the degree of the vertex must be $N-1 \leq d \leq 2\left(N-1\right)$.
If $w=2$ and $d=N-1,N \,\,\text{or}\,\, 2\left(N-1\right)$, then we shall show $x_{d,w} > 0$.
By \eqref{Markinkiewicz}, \eqref{c(n,w)_asz.} and \eqref{A(n,d,w)}, we can compute the asymptotic behaviour of $A \left( n,d,2 \right)$ as follows
$$
A \left( n,N-1,2 \right) \sim p a_{2} \dfrac{n^{2\alpha + \beta +1}}{2\alpha + \beta +1} \left( 1-p \right)  q  x_{N-1,1} \to\infty,
$$
$$
A \left( n,N,2 \right) \sim p a_{2} \dfrac{n^{2\alpha + \beta +1}}{2\alpha + \beta +1} \dfrac{N-1}{N} p r x_{N-1,1} \to\infty,
$$
$$
A \left( n,2\left(N-1\right),2 \right) \sim p a_{2} \dfrac{n^{2\alpha + \beta +1}}{2\alpha + \beta +1} \left[ \left(N-1\right) \left( 1-r \right) + \dfrac{N \left( 1-p \right) \left( 1-q \right) }{p} \right] x_{N-1,1} \to\infty.
$$
Moreover,
$$
B \left( n,d,2 \right) = \O\left( n^{2 \left( \alpha 2+ \beta \right) +1} \right), \text{thus} 
\left(B \left( n,d,2 \right)\right)^{\frac{1}{2}} \log B \left( n,d,2 \right) = \O \left( A \left( n,d,2 \right)\right)
$$
for $d=N-1,N \text{ and } 2\left(N-1\right). $
Therefore in these cases $Z \left( n,d,2 \right) \sim A \left( n,d,2 \right)$ a.s. on $\{A \left( n,d,2 \right) \to \infty\}$ as $n \to \infty$.
It implies that
\begin{equation}
\dfrac{X \left( n,d,2 \right)}{V_n} = \dfrac{Z \left( n,d,2 \right)}{c\left( n,2 \right)V_n} \sim \dfrac{A \left( n,d,2 \right)}{c\left( n,2 \right)V_n} 
\sim \dfrac{p a_{2} \dfrac{n^{2\alpha + \beta +1}}{2\alpha +\beta +1}}{a_{2} n^{2\alpha + \beta }p n} T_{d,2} = \dfrac{T_{d,2}}{2\alpha + \beta +1}
\end{equation}
with appropriate $T_{d,2}$.
Therefore we have
\begin{equation*}
\dfrac{X \left( n,N-1,2 \right)}{V_n} \to \dfrac{\left( 1-p \right) q x_{N-1,1}}{2\alpha + \beta +1} = x_{N-1,2} > 0,
\end{equation*}
\begin{equation*}
\dfrac{X \left( n,N,2 \right)}{V_n} \to  \dfrac{ \dfrac{N-1}{N} p r x_{N-1,1}}{2\alpha + \beta +1} = x_{N,2} > 0,
\end{equation*}
\begin{equation*}
\dfrac{X \left( n,2\left(N-1\right),2 \right)}{V_n}  \to \left[ \left( N-1 \right) \left( 1-r \right) + \dfrac{N \left( 1-p \right) \left( 1-q \right)}{p} \right] \dfrac{x_{N-1,1}}{2\alpha + \beta +1} = x_{2\left(N-1\right),2} > 0,
\end{equation*}
as $n\to\infty$.
So \eqref{X/Vtox} is valid for $w=2$, $d=N-1, N, 2\left(N-1\right)$.

Consider the cases when $N+1 \leq d \leq 2N-3$ and $w=2$.
These cases are different from the previous ones.
By \eqref{A(n,d,w)} and using \textit{Remark 1}, we have
\begin{multline} 
\begin{split}
A & \left( n,d,2 \right) = 
{\E} Z \left( 1,d,2 \right)+\\
&+\sum_{i=2}^{n} c\left( i,w \right)X(i-1,d-m,1) \times \\ 
& \times \left( p \left( 1-r \right)\dfrac {  \binom{d-m}{N-m-1} \binom{V_{i-1}-d+m-1}{m-1}} {\binom{V_{i-1}}{N-1}} \right.
 + \left. \left( 1-p \right)\left( 1-q \right)\dfrac {\binom{d-m}{N-m-1} \binom{ V_{i-1}-d+m-1}{m} } {\binom{V_{i-1}}{N}} \right),
\end{split}
\end{multline}
where $d-m = N-1$.
Using this and the limit of $X(i-1,N-1,1)/V_{i-1}$, we obtain the asymptotic behaviour of $A \left( n,d,2 \right)$ as follows
$$
A \left( n,d,2 \right) \sim \sum_{i=2}^{n} a_{2} i^{2\alpha + \beta} x_{N-1,1} i p 
\left[
\dfrac{\binom{d-m}{N-m-1} \left( N-1\right)! }{\left(m-1\right)! } 
\dfrac{p \left( 1-r \right)}{\left( pi \right)^{N-m}} 
+ 
\dfrac{\binom{d-m}{N-m-1} N! }{m!}
\dfrac{ \left( 1-p \right) \left( 1-q \right)}{ \left( pi \right)^{N-m}}
\right] \sim
$$
$$
\sim a_{2} x_{N-1,1} C \sum_{i=2}^{n} i^{2\alpha + \beta +1 +m -N} \sim
$$
\begin{equation} \label{A(n,5,2)_as}
\sim a_{2} x_{N-1,1} C \dfrac{n^{2\alpha + \beta +2 +m -N}}{2\alpha + \beta +2 +m -N} = \O\left( n^{2\alpha + \beta 
} \right),
\end{equation}
because $N-m\geq 2$. Here $C$ denotes an appropriate constant.
Now we have
\begin{equation*}
\dfrac{X \left( n,d,2 \right)}{V_n} = \dfrac{Z \left( n,d,2 \right)}{c\left( n,2 \right)V_n} = \dfrac{M \left( n,d,2 \right) + A \left( n,d,2 \right)}{c\left( n,2 \right)V_n}.
\end{equation*}
The behaviour of $A \left( n,d,2 \right)$ is given by \eqref{A(n,5,2)_as}.
We denoted by $B \left( n,d,w \right)$ the increasing predictable process in the Doob-Meyer decomposition of $M^2 \left( n,d,w \right)$.
We know that $B \left( n,d,2 \right) = \O\left( n^{ 4\alpha+ 2\beta +1} \right)$ and so 
$ \left(B \left( n,d,2 \right)\right)^{\frac{1}{2}} \log B \left( n,d,2 \right) = \O \left( n^{2\alpha + \beta + \frac{1}{2} + \varepsilon}\right)$
with arbitrary small positive $\varepsilon$.

Applying Proposition VII-2-4 of \cite{neveu}, 
we have
$$
M \left( n,d,2 \right) = \o \left( \left(B \left( n,d,2 \right)\right)^{\frac{1}{2}} \log B \left( n,d,2 \right) \right) = \o \left(  n^{2\alpha + \beta + \frac{1}{2} + \varepsilon} \right)\,\,\,
\text{a.e. on}\,\,\,\{ B \left( n,d,2 \right) \to \infty \}.
$$
Moreover, on the set $\{ B \left( \infty,d,2 \right) < \infty \}$, the sequence $M( n,d,2)$ is a.s. convergent.
So $M( n,d,2) = \o \left(  n^{2\alpha + \beta + \frac{1}{2} + \varepsilon} \right)$ a.s.\,\,.
Therefore, using \eqref{Markinkiewicz} and \eqref{c(n,w)_asz.}, we have
\begin{equation}
\dfrac{X \left( n,d,2 \right)}{V_n} = \dfrac{M \left( n,d,2 \right) + A \left( n,d,2 \right)}{c\left( n,2 \right)V_n} 
\leq C \dfrac{n^{2\alpha + \beta + \frac{1}{2} + \varepsilon}}{n^{2\alpha + \beta}n} = C \dfrac{1}{n^a} \to 0,
\end{equation}
where $n \to \infty$ and $ \dfrac{1}{4} < a < \dfrac{1}{2}$.
So the proposition is valid for $w=1$ and $w=2$.

Suppose that the statement is true for all weights less than $w$ and for all possible degrees. 
First we study the positive limits.
Consider $A(n,d,w)$ in \eqref{A(n,d,w)} and assume that at least one of the coefficients $x_{d,w-1}$, $x_{d-1,w-1}$, $x_{d-\left(N-1\right),w-1}$ is positive.
Then by \eqref{Markinkiewicz}, \eqref{c(n,w)_asz.}, \eqref{A(n,d,w)}  and using the induction hypothesis, we have
\begin{multline*}
\begin{split}
A \left( n,d,w \right) \sim & \sum_{i=2}^{n}  c \left( i,w \right) \left[ X\left(i-1,d,w-1\right) \left( 1-p \right)  q \dfrac{w-1}{i} + \right. \\
&\left. + X\left(i-1,d-1,w-1\right) p r \dfrac{\left(N-1\right)\left( w-1 \right)}{Ni} \right. +\\
& + X\left(i-1,d-\left(N-1\right),w-1\right) \left( \dfrac{p \left( 1-r \right) \left( N-1 \right)}{pi} +
 \left. \dfrac{ \left( 1-p \right) \left( 1-q \right) N}{pi} \right) \right] \sim
 \end{split}
\end{multline*}
\begin{multline*}
\begin{split}
\sim \sum_{i=2}^{n} \left[ c \left( i,w \right) x_{d,w-1} pi \left( 1-p \right)  q \dfrac{w-1}{i} 
+c \left( i,w \right) x_{d-1,w-1} pi p r \dfrac{\left(N-1\right)\left( w-1 \right)}{Ni} \right. +\\
 + c\left( i,w \right) x_{d-\left(N-1\right),w-1} pi \left( \dfrac{ \left( 1-r \right) \left( N-1 \right)}{i} +
 \left. \dfrac{ \left( 1-p \right) \left( 1-q \right) N}{pi} \right) \right] \sim
\end{split}
\end{multline*}
\begin{multline*}
\begin{split}
\sim  \sum_{i=2}^{n} a_{w} i^{\alpha w + \beta} & \left[ x_{d,w-1} p \alpha_{1} \left(w-1\right) 
+ x_{d-1,w-1} p \alpha_{2}   \left( w-1 \right) \right. +\\
& + x_{d-\left(N-1\right),w-1} p\beta \left.\left.
\vphantom{\dfrac{1}{2}} \right)
\vphantom{\dfrac{1}{2}} \right] \sim
\end{split}
\end{multline*}
\begin{equation}  \label{A_as}
\begin{split}
\sim  p a_{w} \dfrac{n^{\alpha w + \beta +1}}{\alpha w + \beta +1} & \left[ \alpha_1 \left(w-1\right) x_{d,w-1} 
+\alpha_2  \left( w-1 \right)x_{d-1,w-1} + \beta x_{d-\left( N-1 \right),w-1}
 \right].
\end{split}
\end{equation}
During the above computation we deleted all terms having asymptotically smaller degree than the others.

First, we examine the case when the limits are positive.
Suppose that there is at least one positive term in \eqref{A_as}. 
Therefore $A(n,d,w) \sim p a_w n^{\alpha w +\beta +1} x_{d,w} \to \infty $ (because $x_{d,w}>0$).
In this case $\left(B \left( n,d,w \right)\right)^{\frac{1}{2}} \log B \left( n,d,w \right) = \O \left( A \left( n,d,w \right)\right)$.
So, using Proposition VII-2-4 of \cite{neveu},  
we have $Z \left( n,d,w \right) \sim A \left( n,d,w \right)$. 
Therefore
\begin{equation}
\dfrac{X \left( n,d,w \right)}{V_n} = \dfrac{Z \left( n,d,w \right)}{c\left( n,w \right)V_n} \sim \dfrac{A \left( n,d,w \right)}{c\left( n,w \right)V_n} 
\sim \dfrac{p a_{w} n^{\alpha w + \beta +1} x_{d,w}}{a_{w} n^{\alpha w + \beta }p n} = x_{d,w} \quad \text{a.s.} \quad \text{as} \quad n \to \infty,
\end{equation}
where, by \eqref{A_as},
\begin{equation*}
x_{d,w} = \dfrac{1}{\alpha w + \beta +1} \left[ \alpha_{1} \left( w-1 \right) x_{d,w-1} + \alpha_{2}\left(  w-1\right)x_{d-1,w-1} +
\beta x_{d-\left(N-1\right),w-1} \right].
\end{equation*}

To handle the case when the limit is $0$, we argue as follows.
Consider the case when the coefficients $x_{d,w-1}$, $x_{d-1,w-1}$, $x_{d-\left(N-1\right),w-1}$ are equal to zero.

By \eqref{A(n,d,w)} and using the induction hypothesis, we have
\begin{multline*}
\begin{split}
A & \left( n,d,w \right) \sim
\sum_{i=2}^{n} \left[ \vphantom{\dfrac{1}{2}} \right. a_{w} i^{\alpha w + \beta} 
\left( \O \left( \dfrac{1}{i^a} \right) +
\sum_{m=2}^{N-2} X(i-1,d-m,w-1) \right. \times \\
& \times \left( p \left( 1-r \right)\dfrac { \binom{d-m}{N-m-1} \left(N-1\right)!} {\left(m-1\right)!} \dfrac{1}{\left(pi\right)^{N-m}} \right.
 + \left. \left. \left( 1-p \right)\left( 1-q \right)\dfrac {\binom{d-m}{N-m-1} N! } {m!} \dfrac{1}{\left(pi\right)^{N-m}} \right) \right) \sim
\end{split}
\end{multline*}
\begin{equation*}
\sim C_1 \sum_{i=2}^{n} a_{w} i^{\alpha w + \beta -a} + 
 C_2 \sum_{i=2}^{n} \sum_{m=2}^{N-2} i^{\alpha w + \beta} x_{d-m,w-1} \dfrac{1}{\left(pi\right)^{N-m-1}}
\end{equation*}
\begin{equation*}
\leq C \dfrac{n^{\alpha w + \beta + 1 - a}}{\alpha w + \beta + 1 - a} + C \dfrac{n^{\alpha w + \beta }}{\alpha w + \beta} =
\end{equation*}
\begin{equation} \label{as_A(n,3w-1,w)}
= \O \left( n^{\alpha w + \beta +1 -a} \right),
\end{equation}
where $C$ denotes an appropriate constant.
So in this case the asymptotic behaviour of $A \left( n,d,w \right)$ is given by \eqref{as_A(n,3w-1,w)}.
On the other hand, $B \left( n,d,w \right) = \O \left( n^{2 \left(\alpha w + \beta \right)+1} \right)$.

Using \eqref{Markinkiewicz} and \eqref{c(n,w)_asz.}, Proposition VII-2-4 of \cite{neveu}
 implies
\begin{equation*}
\dfrac{X \left( n,d,w \right)}{V_n} = \dfrac{Z \left( n,d,w \right)}{c\left( n,w \right)V_n} = \dfrac{M \left( n,d,w \right) + A \left( n,d,w \right)}{c\left( n,w \right)V_n} =
\end{equation*}
\begin{equation}
=\dfrac{\O \left( n^{\alpha w + \beta +1 -a} \right)}{n^{\alpha w + \beta}pn} = \O \left( n^{-a} \right) \to 0 = x_{d,w} \,, \quad {\text{a.s.}}
\end{equation}
So we have obtained the desired result for the case $0$ limit as well.
\hfill$\Box$
\end{pf}
\begin{rem}
We can see that  for each $d$ with $d \geq N$ there exists $w$ such that  $x_{d,w} > 0$.
\end{rem}
\section{The scale-free property for the weights and degrees}
\begin{lem}  \label{xdw}
 Let $p>0$ and define
 $$
 x_{w} = x_{N-1,w} + x_{N,w} + \dots + x_{\left(N-1\right)w,w}
 $$
 for $w=1,2,\dots$\,.
 Then $x_{w}$, $w=1,2,\dots$\,, are positive numbers satisfying the following recurrence:
$$
x_{1} = \dfrac{1}{\alpha + \beta +1},
$$
\begin{equation} \label{rekurziox(w)-re}
x_{w} = \dfrac{\alpha \left( w-1 \right) + \beta}{\alpha w + \beta +1}x_{w-1}, \quad \text{if} \quad w > 1,
\end{equation}
where
\begin{equation*}
\alpha = \left(1-p\right) q + \dfrac{N-1}{N}pr, \quad \beta =  \left( N-1 \right)\left( 1-r \right) + \dfrac{N\left( 1-p \right)\left( 1-q \right)}{p}.
\end{equation*}
$x_{w}$, $w=1,2,\dots$\,, is a discrete probability distribution.
Moreover, $x_{d,w}$, $d=N-1,N,\dots, \left(N-1 \right) w$, $w=1,2,\dots$\,, is a two-dimensional discrete probability distribution.
\end{lem}
\begin{pf}
If $\alpha=0$, then the statement is an obvious consequence of \eqref{rekurzio_x(d,w)}.
Now assume $\alpha\ne 0$.
As $x_{d,w}$ is defined as $x_{d,w}=0$ for $d\notin\{N-1,N,\dots, \left( N-1 \right)w\}$, therefore $x_{w} = \sum_{d} x_{d,w}$.
From the recurrence \eqref{rekurzio_x(d,w)} for $x_{d,w}$, we obtain
\begin{multline*}
x_{w} = \sum_{d=N-1}^{\left( N-1 \right)w} x_{d,w} = \sum_{d} x_{d,w} = \\
 =\dfrac{1}{\alpha w + \beta +1} \left[ \alpha_1 \left( w-1 \right) \sum_{d} x_{d,w-1} + \alpha_2 \left(  w-1\right) \sum_{d} x_{d-1,w-1} +\beta \sum_{d} x_{d-\left(N-1\right),w-1} \right] =
\end{multline*}
$$
 = \dfrac{\alpha \left( w-1 \right) + \beta}{\alpha w + \beta +1}x_{w-1}\,.
$$
Using this recursive formula for $x_w$, we obtain
$$
x_w = x_1 \prod_{j=2}^w \dfrac{\alpha \left( j-1 \right) + \beta}{\alpha j + \beta +1} = 
\dfrac{1}{\alpha + \beta +1} \dfrac{\alpha + \beta}{2 \alpha + \beta +1} \dfrac{2 \alpha + \beta}{3 \alpha + \beta +1} \dots \dfrac{\left( w-1 \right) \alpha + \beta}{w \alpha + \beta +1} =
$$
$$
 = \dfrac{1}{\alpha w + \beta +1} \prod_{j=1}^{w-1} \dfrac{ \frac{\beta}{\alpha} + j}{\frac{\beta +1}{\alpha} + j}
 = \dfrac{1}{\alpha w + \beta +1}   \dfrac{\varGamma \left( w + \frac{\beta}{\alpha} \right)}{\varGamma \left( 1 + \frac{\beta}{\alpha} \right)} 
 \dfrac{\varGamma \left( 1 + \frac{\beta +1}{\alpha} \right)}{\varGamma \left( w + \frac{\beta+1}{\alpha} \right)} =
$$
\begin{equation}  \label{recformx_w}
= \dfrac{\varGamma \left( 1 + \frac{\beta +1}{\alpha} \right)}{\alpha \varGamma \left( 1 + \frac{\beta}{\alpha} \right)}
\dfrac{\varGamma \left( w + \frac{\beta}{\alpha} \right)}{\varGamma \left( w + \frac{\beta+1}{\alpha}+ 1 \right)}.
\end{equation}
Moreover, by \cite{prudnikov}, we have the following formula:
$$
\sum_{k=0}^{n} \dfrac{\varGamma \left( k+a \right)}{\varGamma \left( k+b \right)} = 
\dfrac{1}{a-b+1} \left[ \dfrac{\varGamma \left( n+a+1 \right)}{\varGamma \left( n+b \right)} - \dfrac{\varGamma \left( a \right)}{\varGamma \left( b-1 \right)} \right].
$$ 
Therefore, by some calculation, we obtain $ \sum_{k=1}^{n} x_k \to 1$, as $n \to \infty$.
So $\sum_{w=1}^{\infty} x_w =1$. 
As $\sum_d x_{d,w}=x_w$, so $\sum_{w=1}^{\infty}\sum_{d=N-1}^{\left(N-1\right)w} x_{d,w} =1$ and therefore $x_{d,w}$, $d=N-1,N,\dots, \left(N-1\right)w$, $w=1,2,\dots$\,, is a two-dimensional discrete probability distribution.
\hfill$\Box$
\end{pf}
Let $X \left( n,w \right)$ denote the number of vertices of weight $w$ after $n$ steps.
Next theorem is the scale-free property for the weights.
It is an extension of \textit{Theorem 3.1} in \cite{BaMo1}, see also \textit{Theorem 3.1} of \cite{FIPB}.
\begin{thm}  \label{theorem:scalefreeWeights}
Let $0<p<1$, $q>0$, $r>0$ and $(1-r)(1-q)>0$.
Then for all $w=1,2,\dots$ we have
\begin{equation}
\dfrac{X \left( n,w \right)}{V_n} \rightarrow x_{w}
\end{equation}
almost surely, as $n \rightarrow \infty $, where $x_{w}$, $w=1,2, \dots$\,, are positive numbers satisfying the recurrence \eqref{rekurziox(w)-re}.
Moreover, 
\begin{equation} \label{x_{w}_asz.}
x_{w} \sim C  w^{- \left( 1 + \frac{1}{\alpha} \right)},
\end{equation}
\smallskip
as $w \to \infty$,  with $C =\varGamma \left( 1 + \frac{\beta +1}{\alpha} \right) \big/ \left({\alpha \varGamma \left( 1 + \frac{\beta}{\alpha} \right)}\right) $.
\end{thm}
\begin{pf} 
We have
$$
X \left( n,w \right) = X \left( n,N-1,w \right) + X \left( n,N ,w \right) + \dots + X \left( n,\left(N-1\right)w,w \right).
$$
By Theorem \ref{limX/V},
\begin{equation*}
\dfrac{X \left( n,w \right)}{V_n} \rightarrow x_{w} = x_{N-1,w}+ \dots + x_{\left(N-1\right)w,w}
\end{equation*}
almost surely, as $n \rightarrow \infty $. 
Here each $x_{w}$ is positive.

Using formula \eqref{recformx_w} and the Stirling-formula for the Gamma function, we have
$$
x_w 
= \dfrac{\varGamma \left( 1 + \frac{\beta +1}{\alpha} \right)}{\alpha \varGamma \left( 1 + \frac{\beta}{\alpha} \right)} 
\dfrac{\varGamma \left( w + \frac{\beta}{\alpha} \right)}{\varGamma \left( w + \frac{\beta+1}{\alpha}+ 1 \right)}
 \sim C_0 
 \dfrac{\left(w + \frac{\beta}{\alpha}\right)^{\left(w + \frac{\beta}{\alpha}\right)}}
 {\left(w + \frac{\beta}{\alpha} + \frac{1}{\alpha} + 1\right)^{\left(w + \frac{\beta}{\alpha} + \frac{1}{\alpha} + 1\right)}} =
$$
$$
= C_0 \left( \dfrac{\left(w + \frac{\beta}{\alpha}\right)}{\left(w + \frac{\beta}{\alpha} + \frac{1}{\alpha} + 1\right)}\right)^{\left(w + \frac{\beta}{\alpha}\right)}
\dfrac{1}{\left(w +\frac{\beta}{\alpha} + \frac{1}{\alpha} + 1\right)^{\frac{1}{\alpha} + 1}} 
\sim C w^{- \left( 1 + \frac{1}{\alpha} \right)},
$$
where 
$C_0 =\dfrac{\varGamma \left( 1 + \frac{\beta +1}{\alpha} \right)}{\alpha \varGamma \left( 1 + \frac{\beta}{\alpha} \right)  } \dfrac{1}{\left(\frac{1}{e}\right)^{1 + \frac{1}{\alpha}}} $
and 
$C=\dfrac{\varGamma \left( 1 + \frac{\beta +1}{\alpha} \right)}{\alpha \varGamma \left( 1 + \frac{\beta}{\alpha} \right)  }$.
\hfill$\Box$
\end{pf}
%
Now  we construct a representation of the limiting joint distribution of degrees and weights.

Let $W$ be a random variable with distribution ${\P} \left( W = w\right) = x_w\,\,, w=1,2,\dots $\,\,.
Let $\xi_1 \equiv N-1$ and $\xi_2,\xi_3,\,\dots$ be independent random variables being independent of $W$, too.
For $w \geq 2$ let $\xi_w$ have the following distribution:
\begin{equation*}
{\P} \left( \xi_w = 0\right) = \dfrac{\alpha_1 \left( w-1 \right)}{\alpha \left( w-1 \right) + \beta}, \quad
{\P} \left( \xi_w = 1\right) = \dfrac{\alpha_2 \left( w-1 \right)}{\alpha \left( w-1 \right) + \beta}, \quad
{\P} \left( \xi_w = N-1\right) = \dfrac{\beta}{\alpha \left( w-1 \right) + \beta}\,.
\end{equation*}
Introduce notation $ S_w = \xi_1 + \xi_2 + \dots + \xi_w $\,.

The following representation of the joint distribution of degrees and weights is useful to obtain scale-free property for degrees.
\begin{thm} \label{jointDistr}
${\P} \left(S_W=d,W=w\right) = x_{d,w}$ \,for all\quad $w=1,2,\dots$\,,\,\,\,$d=N-1,N,\dots,\left( N-1 \right)w$.
\end{thm}
\begin{pf} 
If $w=1$ and $d=N-1$ we have
\begin{equation*}
{\P} \left(S_W=N-1,W=1\right) 
 = {\P} \left(\xi_1=N-1,W=1\right) = {\P} \left(W=1\right) = x_1 = x_{N-1,1}\,.
\end{equation*}
If $w=1$ and $d \ne N-1$, then ${\P} \left(S_W=d,W=1\right) = 0 = x_{d,1}$.
\\
If $w=2$ and $d \not \in \{ N-1,N,2\left( N-1 \right) \}$ then we have
\begin{equation*}
{\P} \left(S_W=d,W=2\right) = {\P} \left(S_2=d,W=2\right) = {\P} \left(\xi_1=N-1, \xi_2 = d-\left( N-1 \right), W=2\right) = 0 = x_{d,2}\,.
\end{equation*}
Using the recursion \eqref{rekurziox(w)-re} and the assumption that $\xi_1, \xi_2, \xi_3, \dots$ are independent random variables which are independent of $W$, we have for $w \geq 2$
\begin{equation*}
{\P} \left(S_W=d,W=w\right) = {\P} \left(S_w=d,W=w\right) = {\P} \left(S_w=d \right){\P} \left(W=w\right)=
\end{equation*}
\begin{eqnarray*}
 = \left[{\P} \left(S_{w-1}=d\right){\P}\left( \xi_w=0 \right) + {\P} \left(S_{w-1}=d-1\right){\P}\left( \xi_w=1 \right)+{\P} \left(S_{w-1}=d-\left( N-1 \right) \right){\P} \left(\xi_w=\left( N-1 \right)\right)\right] \times
\end{eqnarray*}
\begin{eqnarray*}
 \times {\P} \left(W=w-1\right)\dfrac{x_w}{x_{w-1}} =
\end{eqnarray*}
\begin{eqnarray*}
 = \left[{\P} \left(S_{w-1}=d\right){\P}\left( W = w-1 \right)\alpha_1 \left(w-1\right) + {\P} \left(S_{w-1}=d-1\right){\P}\left( W = w-1 \right)\alpha_2 \left(w-1\right)+ \right.
 \end{eqnarray*}
\begin{eqnarray*} 
+ \left. {\P} \left(S_{w-1}=d-\left( N-1 \right) \right){\P} \left(W = w-1\right)\beta \right] \dfrac{1}{\alpha w + \beta +1} = 
\end{eqnarray*}
\begin{eqnarray*}
 =  \dfrac{1}{\alpha w + \beta +1} \left[ \alpha_1 \left(w-1\right) {\P} \left(S_W=d,W=w-1\right) + 
  \alpha_2 \left(w-1\right) \left(S_W=d-1,W=w-1\right) +  \right.
 \end{eqnarray*}
 \begin{eqnarray*}
 \left. + \beta {\P} \left(S_W=d-\left( N-1 \right),W=w-1\right)  \right].
 \end{eqnarray*}
Now, we can see that the sequence ${\P} \left(S_W=d,W=w\right)$ satisfies the same recursion \eqref{rekurzio_x(d,w)} as $x_{d,w}$\,. 
\hfill$\Box$
\end{pf}
\begin{thm} \label{Theorem-x(d,w)_as}
 Suppose that $\alpha_1 > 0$ and $\alpha_2 > 0$\,. 
 Then
\begin{equation}  \label{x_dw=}
x_{d,w} = x_w \dfrac{\alpha}{\sqrt{2 \pi \alpha_1 \alpha_2 w}} 
\left[ \exp \left( - \dfrac{\left(d - {\E} S_w\right)^2}{2 {\D}^2 S_w} \right)  + \O \left( w^{-\frac{1}{2}} \right)\right] \,,\,\,\,\text{as}\,\,\,w \to \infty\,,
\end{equation}
 where the error term $\O \left( w^{-\frac{1}{2}} \right)$ does not depend on $d$\,.
\end{thm}
\begin{pf}
We can follow the ideas of the proof of \textit{Theorem 4.2} in \cite{BaMo2}.
Let $w \geq 1$. 
By the definition of the expected value, we have 
\begin{equation*}
{\E} \xi_w = \dfrac{\alpha_2 \left( w-1 \right)}{\alpha \left( w-1 \right) + \beta} + \left( N-1\right) \dfrac{\beta}{\alpha \left( w-1 \right) + \beta}
= \dfrac{\alpha_2}{\alpha} + \dfrac{\left(\left( N-1\right) \alpha -\alpha_2\right) \beta}{\alpha \left(\alpha \left( w-1 \right) + \beta \right)}\,,
\end{equation*}
if $w \geq 2$, hence
\begin{equation*}
{\E} S_w = {\E} \xi_1 + \dots + {\E} \xi_w = w \dfrac{\alpha_2}{\alpha} + \O \left(\log w\right)\,,
\end{equation*}
as $w \to \infty$.
Similarly, by simple computation, we have
\begin{equation} \label{D}
{\D}^2 \xi_w = \dfrac{\alpha_1 \alpha_2}{\alpha^2} + \O \left( \dfrac{1}{w} \right)\,,\,\,\,\,\ \ {\D}^2 S_w = \dfrac{\alpha_1 \alpha_2}{\alpha^2} w  + \O \left( \log w \right)\,,
\end{equation}
as $w \to \infty$.
\\
Now, we can apply Theorem VII.2.5 in \cite{petrov} for $S_w$. 
The conditions  of that theorem are satisfied, therefore we have
\begin{equation}  \label{Petrov}
\sup_{d \in {\Z} } \left| {\D} S_w {\P} \left( S_w = d \right) - \dfrac{1}{\sqrt{2\pi}} \exp \left( - \dfrac{\left( d- {\E} S_w \right)^2}{2 {\D}^2 S_w} \right)\right| = 
\O \left( \dfrac{1}{\sqrt{w}} \right)\,.
\end{equation}
Using \eqref{D} and \eqref{Petrov}, we obtain 
$\left| {\D} S_w - \dfrac{\sqrt{\alpha_1 \alpha_2 w}}{\alpha} \right| {\P} \left( S_w = d \right) = \O \left( w^{-\frac{1}{2}} \right)$.
Therefore, it follows from \eqref{Petrov}, that
\begin{equation} \label{Alk_Petrov}
\sup_{d \in {\Z} } \left| \dfrac{\sqrt{\alpha_1 \alpha_2 w}}{\alpha}\,\, {\P} \left( S_w = d \right) - 
\dfrac{1}{\sqrt{2\pi}} \exp \left( - \dfrac{\left( d- {\E} S_w \right)^2}{2 {\D}^2 S_w} \right)\right| =   \O \left( \dfrac{1}{\sqrt{w}} \right)\,.
\end{equation}
The independence of $W$ and $\xi_i$ implies that
\begin{equation*}
x_{d,w} = {\P}\left( S_W = d, W = w  \right) = {\P} \left( S_w = d \right)x_w\,.
\end{equation*}
Using this in \eqref{Alk_Petrov}, we can obtain the desired result.
\hfill$\Box$ 
\end{pf}
Our last theorem is an extension of \textit{Theorem 4.3} in \cite{BaMo2} (see also \textit{Theorem 3.4} of \cite{FIPB}) to the case of $N$ interactions. 
The theorem shows the scale-free property for the degrees.
\begin{thm}   \label{ThmScaleFreeDegree}
Let $0<p<1$, $q>0$, $r>0$ and $(1-r)(1-q)>0$.
Let us denote by $U\left( n,d \right)$ the number of vertices of degree $d$ after $n$ steps, 
that is $U\left( n,d \right) = \sum_{w : \frac{d}{N-1} \leq w \leq n+1} X \left( n,d,w \right)$\,.
Then, for any $d \geq N-1$ we have
\begin{equation}  \label{u_d}
\dfrac{U\left( n,d \right)}{V_n} \to u_d = \sum_w x_{d,w}
\end{equation}
a.s. as $n \to \infty$\,, where $u_d$, $d = N-1, N, \dots$, are positive numbers. 
Furthermore,
\begin{equation} \label{u_d as}
u_d \sim \dfrac{\varGamma \left( 1 + \frac{\beta + 1}{\alpha} \right)}{\alpha_2 \varGamma 
\left( 1 + \frac{\beta}{\alpha} \right)}\left( \dfrac{\alpha  d}{\alpha_2} \right)^{- \left( 1 + \frac{1}{\alpha} \right)}\,,
\end{equation}
as $d \to \infty$\,.
\end{thm}
\begin{pf} 
By Theorems \ref{limX/V} and  \ref{jointDistr},  $\dfrac{X \left( n,d,w \right)}{V_n}$  converges almost surely to the distribution $x_{d,w}= {\P} \left(S_W=d,W=w\right)$.
But the cardinalities of terms in $\sum_{w : \frac{d}{N-1} \leq w \leq n+1} X \left( n,d,w \right)$ are not bounded when $n \to \infty$. 
However, using that $x_{d,w}$, $d=N-1,N,\dots, \left(N-1 \right)w$, $w=1,2,\dots$ is a proper two-dimensional discrete distribution, 
therefore the convergence of the marginal distributions is a consequence of the convergence of the two-dimensional distributions. 
So we obtain \eqref{u_d}.

To obtain \eqref{u_d as}, we follow the lines of \cite{BaMo2}.
Let
$$
f = \dfrac{\alpha}{\alpha_2} d\,, \quad H =H_d = \left\lbrace w : f-f^{\frac{1}{2}+ \varepsilon} \leq w \leq f+f^{\frac{1}{2}+ \varepsilon}  \right\rbrace\,,
$$
$$
H^{-} =  H_d^{-} = \left\lbrace  w : w < f-f^{\frac{1}{2}+ \varepsilon } \right\rbrace  \,, 
\quad H^{+} = H_d^{+} = \left\lbrace w : w > f+f^{\frac{1}{2}+ \varepsilon}  \right\rbrace\,
$$
with some fixed $0 < \varepsilon < \dfrac{1}{6}$.

Using Hoeffding's exponential inequality (Theorem 2 in \cite{HO}) for $w \in H^{-}$ we have
\begin{equation*}
{\P} \left(S_w = d \right) \leq
{\P} \left(S_w \geq d \right) \leq {\P} \left(S_w - {\E} S_w \geq d - \dfrac{\alpha_2}{\alpha}w - \O \left( \log w \right) \right) \leq
\end{equation*}
\begin{equation*}
\leq \exp \left\lbrace - \dfrac{2}{\left(N-1\right)^2w} \left( d - \dfrac{\alpha_2}{\alpha}w - \O \left( \log w \right) \right)^2  \right\rbrace =
\exp \left\lbrace -\dfrac{2}{\left(N-1\right)^2} \left( \dfrac{\alpha_2}{\alpha} \right)^2  \dfrac{\left( f-w - \O \left( \log w \right) \right)^2}{w}\right\rbrace.
\end{equation*}
Here $w \in H^{-}$ implies that
$$
\left( f-w - \O \left( \log w \right) \right)^2 = \left( f-w \right)^2 - 2 \left( f-w \right) \O \left( \log w \right) + \left(\O \left( \log w \right)\right)^2  \geq f^{1+2\varepsilon} - \O \left( f \log f \right).
$$
Therefore in the case when $w \in H^{-}$ we have
$$
{\P} \left(S_w = d \right) \leq \exp \left\lbrace -\dfrac{2}{\left(N-1\right)^2} \left( \dfrac{\alpha_2}{\alpha} \right)^2  \dfrac{ f^{1+2\varepsilon} -\O \left( f \log f \right)}{f}\right\rbrace =
\exp \left\lbrace - \dfrac{2}{\left(N-1\right)^2}\left( \dfrac{\alpha_2}{\alpha} \right)^2  f^{2\varepsilon} + \O \left(\log f \right) \right\rbrace.
$$
Using this, we can obtain that
$$
{\P} \left(S_W=d,W \in H^{-}\right) = \sum_{w \in H^{-}} {\P} \left(S_w=d,W = w\right) \leq \sum_{w \in H^{-}} {\P} \left(S_w=d\right) \leq
$$
\begin{equation} \label{H-}
\leq f  \exp \left\lbrace - \dfrac{2}{\left(N-1\right)^2}\left( \dfrac{\alpha_2}{\alpha} \right)^2  f^{2\varepsilon} + \O \left(\log f \right) \right\rbrace = 
\o \left( f^{- \left( 1 + \frac{1}{\alpha} \right)} \right).
\end{equation}
Similarly, if $w \in H^{+}$, again by Hoeffding's inequality,  we have
\begin{multline*}
{\P} \left(S_w = d \right) \leq
{\P} \left(S_w \leq d \right) \leq {\P} \left(S_w - {\E} S_w \leq d - \dfrac{\alpha_2}{\alpha}w\right) \leq 
\\ \leq \exp \left\lbrace -\dfrac{2}{\left(N-1\right)^2w} \left(d- \dfrac{\alpha_2}{\alpha}w \right)^2 \right\rbrace =
 \exp \left\lbrace - \dfrac{2}{\left(N-1\right)^2}\left( \dfrac{\alpha_2}{\alpha} \right)^2  \dfrac{\left(f-w\right)^2}{w}\right\rbrace.
\end{multline*}
Using that $w\in H^+$ and $\frac{1}{2} + \varepsilon <1$, we obtain
$2 \left( w-f \right) \geq f^{\frac{1}{2} + \varepsilon} + w - f \geq f^{\frac{1}{2} + \varepsilon} + \left(w-f\right)^{\frac{1}{2} + \varepsilon} \geq w^{\frac{1}{2} + \varepsilon}$ for $d$ large enough. 
Therefore
$$
{\P} \left(S_w = d \right) \leq \exp \left\lbrace -\dfrac{2}{\left(N-1\right)^2} \left( \dfrac{\alpha_2}{\alpha} \right)^2  \dfrac{w^{1+2\varepsilon}}{4w}\right\rbrace = \exp \left\lbrace - \dfrac{1}{2\left(N-1\right)^2}\left( \dfrac{\alpha_2}{\alpha} \right)^2  w^{2\varepsilon} \right\rbrace.
$$
Hence
\begin{equation} \label{H+}
{\P} \left(S_W=d,W \in H^{+}\right) 
\leq \sum_{\{ w \, : \, f < w \}} \exp \left\lbrace - \dfrac{1}{2\left(N-1\right)^2} \left( \dfrac{\alpha_2}{\alpha} \right)^2  w^{2\varepsilon} \right\rbrace  
= \o \left( f^{\left(-1+\frac{1}{\alpha}\right)} \right)
\end{equation}
for $f$ large enough.

Now consider the case when $w \in H$.
First we need some general facts.
Consider the set
$$
B= \left\{ (d,w) \, : \, w\ge 1, d\ge N-1, w\in H_d \right\}.
$$
It is easy to see that when $(d,w)\in B$ then $d\to \infty$ if and only if $w\to \infty$.
More precisely, 
$$
\frac{w}{d} \to 1, \quad {\text{ if}} \quad d\to \infty \quad {\text{and}} \quad  (d,w)\in B.
$$
We have $w =f + \O \left( f^{\frac{1}{2} + \varepsilon}\right)$. 
Then (with $\varepsilon_1 >0$ arbitrarily small) 
\begin{equation}  \label{H++}
- \dfrac{\left(d - {\E} S_w\right)^2}{2 {\D}^2 S_w} =
 - \dfrac{\left(d - w \dfrac{\alpha_2}{\alpha} - \O \left(\log w\right)\right)^2}{2 \dfrac{\alpha_1 \alpha_2}{\alpha^2} w  + \O \left( \log w \right)} =
- \dfrac{\alpha_2}{\alpha_1} \dfrac{\left( f-w  -\O \left(\log w\right) \right)^2}{2w + \O \left( \log w \right)}=
\end{equation}
\begin{equation*}
= - \dfrac{\alpha_2}{\alpha_1} \dfrac{\left( f-w \right)^2 + \O \left( f^{\frac{1}{2} + \varepsilon + \varepsilon_1} \right)}{2w + \O \left( \log w \right)}=
- \dfrac{\alpha_2}{\alpha_1} \dfrac{\left( f-w \right)^2 + \O \left( f^{\frac{1}{2} + \varepsilon + \varepsilon_1} \right)}{2f} 
\dfrac{2f}{2f + \O \left( f^{\frac{1}{2} + \varepsilon} \right)} = 
\end{equation*}
\begin{equation*}
= - \dfrac{\alpha_2}{\alpha_1} 
\dfrac{\left( f-w \right)^2 + \O \left( f^{\frac{1}{2} + \varepsilon + \varepsilon_1} \right)}{2f} 
\left[ 1 - \dfrac{\O \left( f^{\frac{1}{2} + \varepsilon}\right)}{2f + \O \left( f^{\frac{1}{2} + \varepsilon} \right)} \right]
= - \dfrac{\alpha_2}{\alpha_1} \dfrac{\left( f-w \right)^2}{2f} + \O \left( f^{- \frac{1}{2} + 3\varepsilon} \right),
\end{equation*}
as $d \to \infty $.
Here the error term does not depend on $w$.
We shall apply  Theorem \ref{Theorem-x(d,w)_as} that is formula \eqref{x_dw=}.
The asymptotic behaviour of $x_{w}$ is known from \eqref{x_{w}_asz.}. 
Using these facts and \eqref{H++}, we obtain
$$
x_{d,w} \sim C  w^{- \left( 1 + \frac{1}{\alpha} \right)} 
\dfrac{\alpha}{\sqrt{2 \pi \alpha_1 \alpha_2 w}} 
\left[ \exp \left\lbrace - \dfrac{\alpha_2}{\alpha_1} \dfrac{\left( f-w \right)^2}{2f} + \O \left( f^{- \frac{1}{2} + 3\varepsilon} \right) \right\rbrace +
\O \left( w^{-\frac{1}{2}} \right) \right] \sim
$$
$$
\sim C  f^{- \left( 1 + \frac{1}{\alpha} \right)} \dfrac{\alpha}{\alpha_2} \dfrac{1}{\sqrt{2 \pi \frac{\alpha_1}{\alpha_2} f}} \exp \left\lbrace -\dfrac{\left( f-w \right)^2}{2 \frac{\alpha_1}{\alpha_2}f} \right\rbrace
$$
as $d \to \infty$ and $w \in H$, where $C= \varGamma \left( 1 + \frac{\beta + 1}{\alpha} \right)/ \left(\alpha \varGamma \left( 1 + \frac{\beta}{\alpha} \right)\right)$.
Therefore 
$$
\sum_{w \in H} x_{d,w} \sim 
\sum_{f-f^{\frac{1}{2} + \varepsilon} < w < f+f^{\frac{1}{2} + \varepsilon}} C  f^{- \left( 1 + \frac{1}{\alpha} \right)} \dfrac{\alpha}{\alpha_2} 
\dfrac{1}{\sqrt{2 \pi \frac{\alpha_1}{ \alpha_2}f}} \exp \left\lbrace - \dfrac{\left( f-w \right)^2}{2 \frac{\alpha_1}{ \alpha_2}f} \right\rbrace \sim
$$
$$
\sim C  f^{- \left( 1 + \frac{1}{\alpha} \right)} \dfrac{\alpha}{\alpha_2} 
\sum_{-f^{\frac{1}{2}+\varepsilon} < k < f^{\frac{1}{2}+\varepsilon}} \dfrac{1}{\sqrt{2 \pi \frac{\alpha_1}{ \alpha_2}f}} 
\exp \left\lbrace - \dfrac{k^2}{2 \frac{\alpha_1}{ \alpha_2}f} \right\rbrace =
$$
$$
= A \sum_{-f^{\varepsilon} < \frac{k}{\sqrt{f}} < f^{\varepsilon}} \dfrac{1}{\sqrt{f}} \dfrac{1}{\sqrt{2 \pi \frac{\alpha_1}{ \alpha_2}}} 
\exp \left\lbrace - \dfrac{\left(\frac{k}{\sqrt{f}}\right)^2}{2 \frac{\alpha_1}{ \alpha_2}} \right\rbrace 
\to A \int_{-\infty}^{+\infty} \dfrac{1}{\sqrt{2 \pi \frac{\alpha_1}{ \alpha_2}}}  \exp \left\lbrace - \dfrac{x^2}{2 \frac{\alpha_1}{ \alpha_2}} \right\rbrace dx = A.
$$
Thus we have
\begin{equation} \label{H}
{\P} \left(S_W=d,W \in H\right) \sim A =
\dfrac{\varGamma \left( 1 + \frac{\beta + 1}{\alpha} \right)}{\alpha_2 \varGamma \left( 1 + \frac{\beta}{\alpha} \right)}
\left( \dfrac{\alpha  d}{\alpha_2} \right)^{- \left( 1 + \frac{1}{\alpha} \right)},
\end{equation}
as $d \to \infty$.
Finally, from \eqref{H-}, \eqref{H+} and \eqref{H}, we obtain
$$
u_d = \sum_w x_{d,w} = \sum_{w \in H^{-}} x_{d,w} + \sum_{w \in H} x_{d,w} + \sum_{w \in H^{+}} x_{d,w} \sim 
\o \left( f^{- \left( 1+\frac{1}{\alpha}\right)} \right) + C \dfrac{\alpha}{\alpha_2} f^{- \left( 1+\frac{1}{\alpha}\right)} + \o \left( f^{- \left( 1+\frac{1}{\alpha}\right)} \right) \sim
$$
$$
\sim \dfrac{\varGamma \left( 1 + \frac{\beta +1}{\alpha} \right)}{\alpha_2 \varGamma \left( 1 + \frac{\beta}{\alpha} \right)} 
\left(\dfrac{\alpha}{\alpha_2}d\right)^{-\left(1 + \frac{1}{\alpha}\right)},
$$
as $d \to \infty $. 
The proof is complete.
\hfill$\Box$
\end{pf}

\end{document}